\documentclass[a4paper]{article}
\usepackage{amsmath,amssymb,amsfonts}

\title{On the $3 \times 3$ magic square constructed with nine distinct square numbers}

\author{Jailton C. Ferreira}

\date{ }
\begin{document}
\maketitle
\pagenumbering{arabic}

\begin{abstract}
\begin{center}
A proof that there is no  $3 \times 3$ magic square constructed with nine distinct square numbers is given.
\end{center}
\end{abstract}

\section{Introduction} \label{sec-1}

\hspace{22pt} In 1984 Martin Labar \cite{LaBar}  formulated the problem: Can a $3 \times 3$ magic square be constructed with nine distinct square numbers? The problem is found in the second edition of Guy's $\textit{Unsolved}$ $\textit{Problems in Number Theory}$ \cite{Guy}  and became famous when Martin Gardner republished it in 1996 \cite{Gardner}.

\section{The proof}
\label{sec-2}

\begin{center}
\begin{picture}(90, 90)
\linethickness{0.1mm}

\put (0,0) {\line(1,0){90}} \put (0,30) {\line(1,0){90}} \put
(0,60) {\line(1,0){90}} \put (0,90) {\line(1,0){90}}

\put (0,0) {\line(0,1){90}} \put (30,0) {\line(0,1){90}} \put
(60,0) {\line(0,1){90}} \put (90,0) {\line(0,1){90}}

 \put (13,73) {$a$} \put (43,73) {$b$} \put (73,73) {$c$} \put (13,43) {$d$} \put (43,43) {$\varepsilon$} \put (73,43) {$f$} \put (13,13) {$g$} \put (43,13) {$h$} \put (73,13) {$i$} 

\end{picture}
\end{center}
\begin{center}
\textrm{Figure 1}
\end{center}

\hspace{22pt}Let be the square given in Figure 1 such that $a, b, c, d, \varepsilon, f, g, h, i \in \textbf{N}$ and

\begin{equation}\label{dois-1}
a+b+c=x
\end{equation}
\begin{equation}\label{dois-2}
d+\varepsilon+f=x
\end{equation}
\begin{equation}\label{dois-3}
g+h+i=x
\end{equation}
\begin{equation}\label{dois-4}
a+d+g=x
\end{equation}
\begin{equation}\label{dois-5}
b+\varepsilon+h=x
\end{equation}
\begin{equation}\label{dois-6}
c+f+i=x
\end{equation}
\begin{equation}\label{dois-7}
a+\varepsilon+i=x
\end{equation}
\begin{equation}\label{dois-8}
c+\varepsilon+g=x
\end{equation}

The equations \eqref{dois-1}, \eqref{dois-2}, \eqref{dois-3}, \eqref{dois-5}, \eqref{dois-6}  \eqref{dois-4}, \eqref{dois-7} and \eqref{dois-8} can be rewritten as
\begin{equation}\label{dois-9}
a=x-b-c
\end{equation}
\begin{equation}\label{dois-10}
d=x-\varepsilon-f
\end{equation}
\begin{equation}\label{dois-11}
g=x-\varepsilon-i
\end{equation}
\begin{equation}\label{dois-12}
b=x-\varepsilon-h
\end{equation}
\begin{equation}\label{dois-13}
c=x-f-i
\end{equation}

\begin{equation}\label{dois-14}
a + d + g - x=0
\end{equation}
\begin{equation}\label{dois-15}
a + \varepsilon + i - x=0
\end{equation}
\begin{equation}\label{dois-16}
c + \varepsilon + g - x=0
\end{equation}

Substituting sequentially $a$, $d$, $g$, $b$ and $c$ given by  \eqref{dois-10} to \eqref{dois-13} into the equations \eqref{dois-14} to \eqref{dois-16} we obtain, respectively,

\begin{equation}\label{dois-17}
0=0
\end{equation}
\begin{equation}\label{dois-18}
2 \varepsilon + f + h + 2 i - 2 x=0
\end{equation}
\begin{equation}\label{dois-19}
\varepsilon - f - h - 2 i + x=0
\end{equation}

Summing \eqref{dois-18} and \eqref{dois-19} we find
\begin{equation}\label{dois-20}
\varepsilon=\frac{x}{3}
\end{equation}

The set of equations \eqref{dois-1} to \eqref{dois-8} can be put in the form
\begin{equation}\label{dois-21}
a=\varepsilon+\Delta_{1}
\end{equation}
\begin{equation}\label{dois-22}
b=\varepsilon-(\Delta_{1}+\Delta_{2})
\end{equation}
\begin{equation}\label{dois-23}
c=\varepsilon+\Delta_{2}
\end{equation}
\begin{equation}\label{dois-24}
d=\varepsilon-(\Delta_{1}-\Delta_{2})
\end{equation}
\begin{equation}\label{dois-25}
f=\varepsilon+(\Delta_{1}-\Delta_{2})
\end{equation}
\begin{equation}\label{dois-26}
g=\varepsilon-\Delta_{2}
\end{equation}
\begin{equation}\label{dois-27}
h=\varepsilon+(\Delta_{1}+\Delta_{2})
\end{equation}
\begin{equation}\label{dois-28}
i=\varepsilon-\Delta_{1}
\end{equation}
Let us notice that $\Delta_{1} \ne 0$, $\Delta_{2} \ne 0$ and $\Delta_{1} \ne \Delta_{2}$ to obtain the magic square constructed with nine distinct square numbers.

\hspace{22pt} Let us consider that
\begin{equation}\label{dois-29}
\eta_{n}=(n+1)^2-n^2
\end{equation}
\begin{equation}\label{dois-30}
\eta_{n+1}=\eta_{n}+2
\end{equation}
where $n \in \textbf{N}$ and $n>0$. The Figure 2 was obtained using the equations \eqref{dois-29} and \eqref{dois-30}.

\begin{center}
\begin{picture}(200,90)
\linethickness{0.1mm}

\put (15,40) {\line(3,-4){12}} \put (45,40) {\line(3,-4){12}} \put (75,40) {\line(3,-4){12}} \put (105,40) {\line(3,-4){12}} \put (135,40) {\line(3,-4){12}}

\put (45,40) {\line(-3,-4){12}} \put (75,40) {\line(-3,-4){12}} \put (105,40) {\line(-3,-4){12}} \put (135,40) {\line(-3,-4){12}} \put (165,40) {\line(-3,-4){12}}

 \put (13,73) {$1^2$} \put (43,73) {$2^2$} \put (73,73) {$3^2$}  \put (103,73) {$4^2$}  \put (133,73) {$5^2$} \put (163,73) {$6^2$}  \put (193,73) {$...$}

\put (13,43) {$1$} \put (43,43) {$4$} \put (73,43) {$9$}  \put (100,43) {$16$} \put (130,43) {$25$} \put (160,43) {$36$} \put (193,43) {$...$}

\put (28,13) {$3$} \put (58,13) {$5$} \put (88,13) {$7$} \put (118,13) {$9$} \put (145,13) {$11$} \put (193,13) {$...$}

\end{picture}
\textrm{Figure 2}
\end{center}

Let us assume that there exist $0<m<e$ and $e<n$ such that
\begin{equation}\label{dois-31}
n^2+m^2=2 e^2
\end{equation}

The values of $n^2$ and $m^2$ are
\begin{equation}\label{dois-32}
n^2=e^2+\sum_{k=e}^{n-1} \eta_{k}
\end{equation}
and
\begin{equation}\label{dois-33}
m^2=e^2-\sum_{k=m}^{e-1} \eta_{k}
\end{equation}
or
\begin{equation}\label{dois-34}
n^2=e^2+\frac{\eta_{e}+\eta_{n-1}}{2}((n-1)-(e-1))
\end{equation}
and
\begin{equation}\label{dois-35}
m^2=e^2-\frac{\eta_{m}+\eta_{e-1}}{2}((e-1)-(m-1))
\end{equation}

Considering $a=n^2$, $i=m^2$ and $\varepsilon=e^2$ we have
\begin{equation}\label{dois-36}
\frac{\eta_{e}+\eta_{n-1}}{2}((n-1)-(e-1))-\frac{\eta_{m}+\eta_{e-1}}{2}((e-1)-(m-1))=0
\end{equation}
or
\begin{eqnarray}\label{dois-37}
\lefteqn{(-e+n)(-2e^2+2(1+e)^2+2(-1-e+n))-{} }
\nonumber\\
& & {}((-2-2e^2+2(1+e)^2-2(e-m))(e-m)=0
\end{eqnarray}

Let us assume that there exist $w$ and $z$ such that
\begin{equation}\label{dois-38}
c=(n+w)^2
\end{equation}
and
\begin{equation}\label{dois-39}
g=(m-z)^2
\end{equation}
where $w$ and $z$ are positive integers. In this case we have
\begin{eqnarray}\label{dois-40}
\lefteqn{(-e+n+w)(-2e^2+2(1+e)^2+2(-1-e+n+w))-{} }
\nonumber\\
& & {}((-2-2e^2+2(1+e)^2-2(e-m+z))(e-m+z)=0
\end{eqnarray}

Subtracting \eqref{dois-37} from \eqref{dois-40} we obtain
\begin{equation}\label{dois-41}
2n w+w^2+(-2m+z)=0
\end{equation}
Solving \eqref{dois-41} for $z$ we find
\begin{equation}\label{dois-42}
z_{1}=m+\sqrt{m^2-2n w-w^2}
\end{equation}
and
\begin{equation}\label{dois-43}
z_{2}=m-\sqrt{m^2-2n w-w^2}
\end{equation}
The root $z_{1}$ implies that $m-z$ is not a positive integer, however $m-z$ must be a positive integer. Therefore $z=z_{2}$.

\hspace{22pt} We have assumed
\begin{equation}\label{dois-44}
 2 e^2=m^2+n^2
\end{equation}
and
\begin{equation}\label{dois-45}
 2 e^2=(m-z)^2+(n+w)^2
\end{equation}
Subtracting \eqref{dois-44} from \eqref{dois-45} we obtain
\begin{equation}\label{dois-46}
(m-z)^2+(n+w)^2 -(m^2+n^2)=0
\end{equation}
Substituting \eqref{dois-43} into \eqref{dois-46} we have
\begin{equation}\label{dois-47}
n^2-2n w-w^2-(n+w)^2=0
\end{equation}

Solving \eqref{dois-47} for $w$ we find
\begin{equation}\label{dois-48}
w_{1}=0
\end{equation}
and
\begin{equation}\label{dois-49}
w_{2}=-2n
\end{equation}

The root $w_{2}$ implies in $n+w$ negative, however $n+w$ must be a positive integer. Therefore $w=0$. Since that $w=0$ the condition $\Delta_{1} \ne \Delta_{2}$ is not satisfied.There is no magic square constructed with nine distinct square numbers.

\end{document}